\documentclass[11pt]{article}
\newcommand{\version}{November 30, 2007 }

\newcommand{\lanbox}{\hfill \hbox{$\, 
\vrule height 0.25cm width 0.25cm depth 0.01cm
\,$}}

\usepackage{amsthm,amsfonts,amsmath, amscd}
\swapnumbers
\theoremstyle{plain}

\newtheorem{thm}{THEOREM}[section]
\newtheorem{lm}[thm]{LEMMA}

\theoremstyle{definition}

\theoremstyle{definition}
\newtheorem{remark}[thm]{Remark}
\newcommand{\upchi}{\raise1pt\hbox{$\chi$}}
\newcommand{\R}{{\mathord{\mathbb R}}}

\newcommand{\C}{{\mathord{\mathbb C}}}

\newcommand{\hn}{{\mathord{\widehat{n}}}}

\newcommand{\tr}{{\rm Tr}}
\renewcommand{\|}{{\Vert}}
\numberwithin{equation}{section}
\begin{document}

\markboth{\scriptsize{CL \version}}{\scriptsize{CL \version}}

\def\mn{{\bf M}_n}
\def\hn{{\bf H}_n}
\def\hnp{{\bf H}_n^+}
\def\hmnp{{\bf H}_{mn}^+}
\def\h{{\cal H}}

\title{A MINKOWSKI TYPE TRACE INEQUALITY AND STRONG SUBADDITIVITY OF QUANTUM ENTROPY II:\\ CONVEXITY AND CONCAVITY}
\author{\vspace{5pt} Eric A. Carlen$^1$ and
Elliott H. Lieb$^{2}$ \\
\vspace{5pt}\small{$1.$ Department of Mathematics, Hill Center,}\\[-6pt]
\small{Rutgers University,
110 Frelinghuysen Road
Piscataway NJ 08854-8019 USA}\\
\vspace{5pt}\small{$2.$ Departments of Mathematics and Physics, Jadwin
Hall,} \\[-6pt]
\small{Princeton University, P.~O.~Box 708, Princeton, NJ
  08544}\\
 }
\date{\version}
\maketitle
\footnotetext                                                                         
[1]{Work partially
supported by U.S. National Science Foundation
grant DMS 06-00037.    }                                                          
\footnotetext
[2]{Work partially
supported by U.S. National Science Foundation
grant PHY 06-52854.\\
\copyright\, 2007 by the authors. This paper may be reproduced, in its
entirety, for non-commercial purposes.}

\begin{abstract}

We revisit and prove some convexity inequalities  
for trace functions conjectured
 in the earlier part I.  The main functional considered is 
$$
\Phi_{p,q} (A_1,\, A_2,\ \dots, A_m) = \left(\tr\left[( \, {\textstyle\sum_{j=1}^m
 A_j^p } \, )
^{q/p} \right] \right)^{1/q}
$$
for $m$ positive definite operators $A_j$. In part I we only
considered the case $q=1$ and proved the concavity of $\Phi_{p,1}$ for
$0 <p\leq 1$ and the convexity for $p=2$. We conjectured the convexity
of $\Phi_{p,1}$ for $1< p < 2$.  Here we not only settle the
unresolved case of joint convexity for $1\leq p \leq 2$, we are also
able to include the parameter $q\geq 1$ and still retain the
convexity.  Among other things this leads to a definition of an
$L^q(L^p)$ norm for operators when $1\leq p \leq 2$
and a Minkowski inequality for operators on a tensor product of 
three  Hilbert spaces -- which leads to another proof of
strong subadditivity of entropy. We also prove
convexity/concavity properties of some other, related functionals. 

\end{abstract}

\medskip
\leftline{\footnotesize{\qquad Mathematics subject classification numbers:  47A63, 15A90}}
\leftline{\footnotesize{\qquad Key Words: convexity, concavity, trace inequality, 
entropy, operator norms}}

\section{Introduction} \label{intro}

Let $\mn$ denote the set of $n\times n$ matrices  with complex
entries, and let $A^*$ denote the Hermitian conjugate of $A\in \mn$.
For $0< p < \infty$, and $A\in \mn$, define
$$
\|A\|_q = ({\rm Tr}[(A^*A)^{q/2}])^{1/q}\ .
$$
For $q \ge 1$, this defines a norm on $\mn$, but not for $q<1$.
Nonetheless, it will be convenient here to use this notation for all
$q>0$.  Finally, define $\|A\|_\infty$ to be the operator norm of $A$.

Let $\hn$ denote the set of $n\times n$ Hermitian matrices, and let
$\hnp$ denote the set of $n\times n$ positive semidefinite matrices.  For any two
positive integers $m$ and $n$, we may identify operators on
$\C^m\otimes \C^n$ with $mn\times mn$ matrices, and then we have the
two partial traces ${\rm Tr}_1$ and ${\rm Tr_2}$: For any $mn\times
mn$ matrix $A$, ${\rm Tr}_2A$ is the $m\times m$ matrix such that for
all $m\times m$ matrices $B$,
$$
{\rm Tr}[({\rm Tr}_2A)B]  = {\rm Tr}[A(B\otimes I_{n\times n})]\ ,
$$
and ${\rm Tr}_1$, which ``traces out'' the first factor, is defined
analogously.  Later, this notation will be extended to tensor products
of more factors in the obvious way, without further discussion.

This paper concerns  properties of the following trace functionals:

\medskip
\noindent{\it (1)} For any  numbers $p,q>0$, and any positive integer $m$, define  
$\Phi_{p,q}$ on the $m$--fold Cartesian product of $\hnp$ with itself by

\begin{equation}
\boxed{\Phi_{p,q}(A_1,\dots,A_m) 
= \|\  ({\textstyle
\sum_{j=1}^m A_j^p
)^{1/p}  } \ \|_q\ .}
\end{equation}

\medskip
\noindent{\it (2)} For any  numbers $p,q>0$, and any positive integers $m$  and $n$, define  
$\Psi_{p,q}$  on $\hmnp$  by  
\begin{equation}
\boxed{\Psi_{p,q}(A) = \| \left({\rm Tr}_2A^p\right)^{1/p}\|_q\ .}
\end{equation}

\medskip
\noindent{\it (3)} For any fixed $n\times n$ matrix $B$ and any numbers $p,q>0$, define
$\Upsilon_{p,q}$ on $\hnp$ by
\begin{equation}
\boxed{\Upsilon_{p,q}(A) = {\rm Tr}\left[(B^*A^pB)^{q/p}\right]\ .}
\end{equation}

\medskip

We shall determine conditions on $p$ and $q$ under which these
functionals are convex or concave (jointly, in the case of
$\Phi_{p,q}$).

\medskip
\begin{thm}\label{main} 
For all $1  \le p \le 2$,  and for all $ q\geq 1$, $\Upsilon_{p,q}$  and $\Psi_{p,q}$ are convex on $\hnp$
and $\hmnp$ respectively, while $\Phi_{p,q}$ is jointly convex on $(\hnp)^m$.

For $0 \le p \le q \le 1$,  $\Upsilon_{p,q}$  and $\Psi_{p,q}$ are concave on $\hnp$
and $\hmnp$ respectively, while $\Phi_{p,q}$ is jointly concave on $(\hnp)^m$.

For $p>2$, none of these functions are convex or concave for any values of $q\neq p$. 

\end{thm}

\remark \label{r1}
Note that $\Upsilon_{p,q}$ is homogeneous of degree $q\geq 1$. By a 
general argument  $\left(\Upsilon_{p,q}\right)^{1/q}$ 
is also convex/concave.  (The general argument is this: 
A function $f$ that is homogeneous
of degree one is convex (concave) if and only if the level set $\{x\
:\ f(x) \le 1\}$ ( $\{x\ :\ f(x) \ge 1\}$) is convex. Hence,  if $g(x)$
is homogeneous of degree $q$, and convex (concave), so that $\{x\ :\
g(x) \le 1\}$ ($\{x\ :\ f(x) \ge 1\}$) is convex, $g^{1/q}$ is convex
(concave).)  

\medskip
Apart from the elementary case $p=2, q=1$ which is discussed and proved  in
\cite{CL}, the parts of this theorem that refer to convexity are new.
The parts of this theorem that refer to concavity are already known in
the case $q=1$, but are new for other values of $q$: For $0\le p \le
1$, the concavity of $\Upsilon_{p,1}$ is a theorem of Epstein
\cite{E}, while the concavity of $\Psi_{p,1}$ and joint concavity of
$\Phi_{p,1}$ was proved in \cite{CL}. Also in \cite{CL}, the convexity
of $\Phi_{p,1}$ and $\Psi_{p,1}$ for $1 \le p < 2$ was conjectured.
In Theorem \ref{main}, this conjecture, and more, is verified.

We also note that Bekjan \cite{B} has proved the concavity of
$\Phi_{p,1}$ for $-1 < p < 0$. While our interest here is in positive
values of $p$, we shall return to Bekjan's work towards the end of the
paper when we summarize the relations among these various convexity
and concavity results.

We have stated these inequalities for matrices, but as none of them
refers to the dimension, it is easy to extend them to operators on a
separable Hilbert space $\h$, and we take this for granted in the next
paragraphs.  There are other,  much more interesting extensions to be
considered: For example, in \cite{H}, among other things Hiai has
extended the concavity results from \cite{CL} to a von Neuman algebra
setting in which the trace is replaced by a more general state.  In
this regard see also the paper \cite{K2} by Kosaki.


\subsection {Application to non-commutative Minkowski inequalities}
As an application of Theorem \ref{main}, we deduce a tracial analog of
Minkowski's inequality for multiple integrals: Let $(X,{\rm d}\mu)$
and $(Y,{\rm d}\nu)$ be sigma-finite measure spaces, and then for any
non-negative measurable function $f$ on $X\times Y$, and any $1 \le q
\le p$,
$$
\left[\int_Y\left(\int_X f^q(x,y){\rm d}\mu\right)^{p/q}{\rm
    d}\nu\right]^{1/p}\leq 
\left[\int_X\left(\int_Y f^p(x,y){\rm d}\nu\right)^{q/p}{\rm d}\mu\right]^{1/q} 
\ . 
$$ 
For $0 \le p \le q \le 1$, the inequality
reverses.  (When $p=1$, one has the ``standard'' Minkowski inequality.
This form involving $p$ and $q$ is obtained by applying the standard
inequality  for the $L^{p/q}$ norm to $F(x,y) = f^q(x,y)$.)

Introducing a third measure space $(Z,{\rm d}\lambda)$, we have the
pointwise inequality
\begin{equation}\label{mink31}
\left[\int_Y\left(\int_X f^q(x,y,z){\rm d}\mu\right)^{p/q}{\rm d}\nu
\right]^{1/p}\leq
\left[\int_X\left(\int_Y f^p(x,y,z){\rm d}\nu\right)^{q/p}{\rm d}\mu\right]^{1/q} 
\end{equation}
for any non-negative measurable function $f$ on $X\times Y\times Z$,
and any $1 \le q \le p$.  
Now,  raising both sides to the power $q$ and then integrating over $Z$, one obtains
\begin{equation}\label{mink32}
\int_Z \left[\int_Y\left(\int_X f^q(x,y,z){\rm d}
\mu\right)^{p/q}{\rm d}\nu\right]^{q/p}{\rm d}\lambda \  \leq 
\int_Z \int_X\left(\int_Y f^p(x,y,z)
{\rm d}\nu\right)^{q/p}{\rm d}\mu {\rm d}\lambda  \ .
\end{equation}

To formulate the corresponding tracial inequalities, let ${\cal H}_j$,
$j=1,2,3$ be three separable Hilbert spaces.  Let ${\cal H}$ denote
${\cal H}_1\otimes {\cal H}_2\otimes {\cal H}_3$, and let ${\rm Tr}_j$
be the partial trace that ``traces out'' the $j$th Hilbert space.

Then the tracial analog of (\ref{mink31}) would be
\begin{equation} \label{mink2}
\left({\rm Tr}_2\left[ \left({\rm Tr}_1 A^q\right)^{p/q}\right] \right)^{1/p} \leq
\left({\rm Tr}_1\left[ \left({\rm Tr}_2A^p\right)^{q/p}\right] \right)^{1/q} 
\end{equation}
as an operator inequality on ${\cal H}_3$. When the dimension of
$\h_3$ equals 1, there are really only two spaces present (1 and 2), and then
(\ref{mink2}) is the two space operator inequality proved in \cite{CL}. 
When ${\rm dim}({\cal H}_3)>1$, (\ref{mink2}) cannot hold as an
operator inequality, for any $p>1$, even for $q=1$, as explained in
\cite{CL}.  However, for $1\le p \le 2$, the analog of (\ref{mink32})
is true no matter what ${\rm dim}({\cal H}_3)$ may be:

\medskip
\begin{thm}\label{mink} 
For $1 \le q \le p \le 2$,  and all positive operators $A$ on ${\cal H} = 
{\cal H}_1\otimes {\cal H}_2\otimes {\cal H}_3$,
\begin{equation}\label{minkpq}
{\rm Tr}_3\left({\rm Tr}_2 \left[\left({\rm Tr}_1 A^q\right)^{p/q}\right]\right)^{q/p} \leq
{\rm Tr}_3\left({\rm Tr}_1 \left[\left( {\rm Tr}_2 A^p\right)^{q/p} \right]\right)
\end{equation}
For $0 \le p \leq 1 $, and any $q  \geq p$,   this inequality reverses.
\end{thm}

\medskip For $0 < p < 1$, $q=1$, and for $p=2$, $q=1$, this was proved
in \cite{CL}, while for $1< p < 2$ and $q=1$ it was conjectured in
\cite{CL}.

Note that for $q=p=1$, we have an identity, and so we may obtain a new
inequality by differentiating in $p$ and $q$ at $p=q=1$. As shown in
\cite{CL}, one obtains the strong subadditivity of the
quantum entropy in this way. To state this result, first proved in \cite{LR}, we
recall the definition of the quantum entropy: A {\it density matrix}
$\rho$ on a separable Hilbert space ${\cal H}$ is a positive operator
on ${\cal H}$ such that ${\rm Tr}(\rho) =1$. The {\it entropy} of
$\rho$, $S(\rho)$, is defined by
$$S(\rho) = -{\rm Tr}\left(\rho\ln\rho \right)\ .$$

Since (\ref{minkpq}) is an equality at $p=q =1$, we can obtain
inequalities by differentiating in $p$ at $p=q=1$.   Since
$$\frac {{\rm d}} { {\rm d}p} {\rm Tr}(\rho^p)\big|_{p=1} =  {\rm Tr}(\rho \ln \rho)\ ,$$
the resulting inequality will involve the entropies of various partial
traces of $\rho$. We shall use the following notation for these:
$$
\rho_{123} = \rho \,\quad \rho_{23} = {\rm Tr}_1\rho\ ,\quad \rho _3 
= {\rm Tr}_1{\rm Tr}_2 \rho
$$
and so forth. (That is, the subscripts indicate the spaces ``remaining'' 
after the traces.)

Then differentiating in $p$ at  $p=q =1$, we obtain
the strong subadditivity of  the quantum entropy \cite{LR}, \cite{L}:
$$
S(\rho_{13}) + S(\rho_{23})  \ge S(\rho_{123}) + S(\rho_3)\ ,
$$
just as in \cite{CL}, except that there, only a left derivative was
taken since the Minkowski inequality was then known only for for $p<1$ and
$q=1$.

\remark \label{rssa}  
Since we did not provide the details of the differentiation 
argument in our previous paper, we take
the opportunity to do so here. The basic fact is that for a 
positive operator $A$, and $\varepsilon$
close to zero,
$$
A^{1+\varepsilon} = A + \varepsilon A\ln A + {\cal O}(\varepsilon^2)\ .
$$
At least in finite dimensions, one can take a partial trace of both sides, and the
resulting identity still holds.

Applying this with $A = \rho$ in finite dimensions, we compute
$$
{\rm \tr}_2\left(\rho^{1+\varepsilon}\right) = \rho_{13} + \varepsilon \tr_2(\rho\ln \rho) + {\cal O}(\varepsilon^2)\ .
$$
Then, since to leading order in $\varepsilon$, $1/(1+\varepsilon)$ is $1 - \varepsilon$,
$$
\left[{\rm \tr}_2(\rho^{1+\varepsilon}) \right]^{1/(1+\varepsilon)}
 = \rho_{13} + \varepsilon \tr_2(\rho\ln \rho)  - \varepsilon 
\rho_{13}\ln\rho_{13} +
{\cal O}(\varepsilon^2)\ .
$$
Thus,
\begin{equation}\label{ssa1}
\tr_{1} \left(\tr_{3}\left[{\rm \tr}_2(\rho^{1+\varepsilon})\right]^{1/(1+\varepsilon)} \right)
 = 1 - \varepsilon S(\rho) + \varepsilon S(\rho_{13}) + {\cal O}(\varepsilon^2)\ .
 \end{equation}

In the same way, we find
$$
\tr_2 \left(\rho_{23}^{1+\varepsilon}\right) = \rho_3 + 
\varepsilon\tr_2\left(\rho_{23}\ln \rho_{23}\right) 
 + {\cal O}(\varepsilon^2)\ ,
$$
and then

\begin{equation}\label{ssa2}
\tr_3\left[\tr_2 \left(\rho_{23}^{1+\varepsilon}\right)^{1/(1+\varepsilon)}
\right] = 
1 - \varepsilon S(\rho_{23}) + \varepsilon S(\rho_3) + {\cal O}(\varepsilon^2)\ .
 \end{equation}

Combining (\ref{ssa1}), (\ref{ssa2}) and the $q=1$ case of Theorem 
\ref{mink}, we obtain
$$
S(\rho_{13}) +  S(\rho_{23}) \ge S(\rho) + S(\rho_3)\ 
$$
in the finite dimensional case, where there is no issue of uniformity  in the
${\cal O}(\varepsilon^2)$ remainders; 
The general case follows by finite dimensional approximation (see the appendix
to \cite{LR}.) \lanbox


\subsection{Non-commutative $\mathbf{ L^q(L^p)}$ norms}

We  point out that the convexity 
of $\Psi_{p,q}$ can be used to define certain non-commutative analogs
of the $L^q(L^p) $  norms: Let $(X,{\rm d}\mu)$ and $(Y,{\rm d}\nu)$ be
sigma-finite measure spaces. Any non-negative measurable
function $f$ on $X\times Y$ such that
\begin{equation}  \label{pqnorm}
\left[\int_X\left(\int_Y f^p(x,y){\rm d}\nu\right)^{q/p}
{\rm d}\mu\right]^{1/q} < \infty
\end{equation}
may be regarded as a function on $X$ with values in $L^p(Y,{\rm
  d}\nu)$, and then the quantity above defines the natural $L^q$ norm
of this $L^p$-valued function.

Similarly, a non-negative operator $A$ on the tensor product ${\cal
  H}_1\otimes{\cal H}_2$ of two separable Hilbert spaces can be
regarded as a non-commutative analog of such a function. The problem of
constructing an  $L^q(L^p)$ norm was actively considered by  several 
people in the operator space community. Pisier successfully found one solution
to the problem, valid for all $1\leq q\leq p $. With  $1/r := 1/q -1/p$ 
Pisier's definition is 
\begin{equation} \label{pisier}
\|X\|_{L^q(L^p)} = \inf_{X =
(A\otimes I)Y(B\otimes I)} \|A\|_{2r}\|Y\|_{q}\|B\|_{2r}  \  .
\end{equation}
Pisier proved this is a norm in \cite{P1}. However,  
the definition does not bear such an obvious 
similarity to (\ref{pqnorm}), and it was an actively considered problem
to decide if $\Psi_{p,q}$ is convex and whether this convexity could be used to define 
a  non-commutative    $L^q(L^p)$  norm. The convexity is now established in 
Theorem \ref{main} but, since $\Psi_{p,q}$   is defined only for
positive operators, it does not directly define a norm -- even for self-adjoint operators.
(Note that we cannot simply replace $X$ by $  |X|$ since $X \mapsto  |X|$  is not 
operator convex.)  However, the following standard construction does give us a norm
for {\it self-adjoint} operators when $1\leq p \leq 2$ and all $q\ge 1$:
\begin{equation} \label{psinorm}
|\!|\!| X |\!|\!| _{L^q(L^p)} = \inf_{\scriptsize{\begin{array}{c} X=A-B, \\ A\geq 0, B\geq 0\\ \end{array}} }\left(
\Psi_{p,q}(A)    +\Psi_{p,q}(B) \right) \ .
\end{equation}
This definition satisfies the triangle inequality because if we write
$X=A-B$ and $Y=C-D$ for $A,B,C,D\ge 0$, then $X+Y = (A+C)-(B+D)$ and
\begin{eqnarray}
 |\!|\!| X +Y|\!|\!| _{L^q(L^p)} 
=\inf \cdots &\leq & \Psi_{p,q}(A+C)    +\Psi_{p,q}(B+D)  \nonumber \\ 
&\leq &
\Psi_{p,q}(A)    +\Psi_{p,q}(C) +\Psi_{p,q}(B)    +\Psi_{p,q}(D) \nonumber 
\end{eqnarray} 
The second inequality above is the convexity of $\Psi_{p,q}$, which, since 
$\Psi_{p,q}$ is homogeneous of degree one, amounts to subadditivity.

If we now choose $A,B,C$ and $D$ so that
$\Psi_{p,q}(A)+ \Psi_{p,q}(B)$ nearly equals $ |\!|\!| X |\!|\!| _{L^q(L^p)}$ and 
$\Psi_{p,q}(C)+ \Psi_{p,q}(D)$ nearly equals $ |\!|\!| Y |\!|\!| _{L^q(L^p)}$, we conclude
$$|\!|\!| X +Y|\!|\!| _{L^q(L^p)}  \le |\!|\!| X|\!|\!| _{L^q(L^p)}  + |\!|\!|Y|\!|\!| _{L^q(L^p)} \ .$$
\smallskip

\remark \label{r2}
 It is to be noted that $A \to \Psi_{p,q}(A)$ is {\it not
monotone}, but the norm $ A \to |\!|\!| A |\!|\!|$ {\it is monotone}
(by construction).  Even for $A>0 $, it is not necessarily the case that 
$  |\!|\!| A |\!|\!| = \Psi_{p,q}(A)$ although  
$  |\!|\!| A |\!|\!| \leq \Psi_{p,q}(A)$, in general. 

\medskip

Now that one has an $L^q(L^p)$ norm defined on self adjoint operators $A$ on ${\cal H}_1\otimes{\cal H}_2$,
 one can go on to define an $L^q(L^p)$ norm for {\it all} operators $A$ on ${\cal H}_1\otimes{\cal H}_2$,
self adjoint or not, as follows:  
Consider the block matrix
\begin{equation}\label{block}
{\cal A} = \left[\begin{array}{ccc} 0 & A \\ A^* & 0 \\ \end{array}\right]\ ,
\end{equation}
 which is a self adjoint operator on
${\cal H}_1\otimes{\cal H}_2\otimes \C^2$. For present purposes, we regard this Hilbert space
as ${\cal H}_1\otimes ({\cal H}_2\otimes \C^2)$. That is, when computing $\Psi_{p,q}$ of a positive  
operator on  ${\cal H}_1\otimes ({\cal H}_2\otimes \C^2)$, ${\cal H}_1$ is the first factor, and 
${\cal H}_2\otimes \C^2$ is the second factor.  We then define
\begin{equation}
\boxed{ \|A\|_{L^q(L^p)}  = \frac{1}{2}|\!|\!| {\cal A} |\!|\!|
_{L^q(L^p)}\ ,}
\end{equation}
with $A$ and ${\cal A}$ related as in (\ref{block})

Notice that if $A$ happens to be  self adjoint, and $A = B-C$ where $B$ and $C$ are positive, then 
$${\cal A} = \left[\begin{array}{ccc} 0 & B-C \\ B-C & 0 \\ \end{array}\right] =  \frac{1}{2}\left[\begin{array}{ccc} B+C & B-C \\ B-C & B+C \\ \end{array}\right]
-  \frac{1}{2}\left[\begin{array}{ccc} B+C & C-B \\ C-B & B+C \\ \end{array}\right]$$
is a decomposition of ${\cal A}$ as a difference of positive operators. That is, every decomposition of $A$ as a difference of positive operators induces a corresponding decomposition of ${\cal A}$.

Moreover, with the  unitary operator ${\cal U}$ defined by
${{\cal U} = \displaystyle \frac{1}{\sqrt{2}}\left[\begin{array}{ccc} I & \phantom{-}I \\ I & -I \\ \end{array}\right] }$, 
$$ {\cal U}  \left[\begin{array}{ccc} B+C & B-C \\ B-C & B+C \\ \end{array}\right]{\cal U} ^* = 
2\left[\begin{array}{ccc} B & 0 \\ 0 & C \\ \end{array}\right] \quad{\rm and}\quad
{\cal U}  \left[\begin{array}{ccc} B+C & C-B \\ C-B & B+C \\ \end{array}\right]{\cal U} ^* = 
2\left[\begin{array}{ccc} C & 0 \\ 0 & B \\ \end{array}\right] \ .$$
Because  the unitary operator ${\cal U}$ respects the particular structure of  ${\cal H}_1\otimes ({\cal H}_2\otimes \C^2)$,
$$\Psi_{p,q}\left({\cal U}  \left[\begin{array}{ccc} B+C & B-C \\ B-C & B+C \\ \end{array}\right]{\cal U} ^*\right) = 
2\Psi_{p,q}\left(\left[\begin{array}{ccc} B & 0 \\ 0 & C \\ \end{array}\right] \right) = 2\left(\Psi_{p,q}(B) + \Psi_{p,q}(C)\right)\ ,$$
and similarly for the other term. Thus, 
$$\frac{1}{2}\left( \Psi_{p,q}\left(\left[\begin{array}{ccc} B+C & B-C \\ B-C & B+C \\ \end{array}\right] \right) + 
\Psi_{p,q}\left(\left[\begin{array}{ccc} B+C & C-B \\ C-B & B+C \\ \end{array}\right] \right)\right) = 
2\left(\Psi_{p,q}(B) + \Psi_{p,q}(C)\right)\ .$$
It therefore follows, for self-adjoint $A$,  that
$$\|A\|_{L^q(L^p)}  \le |\!|\!| A |\!|\!| _{L^q(L^p)}\ .$$
In fact, there  are many more ways to decompose ${\cal A}$ as a difference of two positive operators besides just those induced by a decomposition of $A$. 
Therefore, in general we can expect the last inequality to be strict. 

We note that we could have proceeded slightly differently: For the purpose of computing $\Psi_{p,q}$ of a positive operator on 
${\cal H}_1\otimes{\cal H}_2\otimes \C^2$, we could have regarded this Hilbert space as  $(\C^2\otimes {\cal H}_1)\otimes {\cal H}_2$. 
This provides a second, distinct, way to construct the norm for general operators. It is an interesting open problem to compare the norms
constructed this way using $\Psi_{p,q}$ with Pisier's  non-commutative $L^q(L^p)$ norm defined in (\ref{pisier}).

\subsection{Organization of the paper} 
In the next sections, we give
the proof of Theorem \ref{main}. The key will be to prove the
convexity and concavity properties of $\Upsilon_{p,q}$, which amounts
to a generalization of Epstein's Theorem.  This will be done in
Section 2.

Then,  in Section 3, we shall be able to deduce the convexity and
concavity properties of the other two functionals from this using the
same strategy that was used to deduce the corresponding concavity results
from Epstein's Theorem in \cite{CL}.  The proof of Theorem \ref{main}
concludes Section 3.  In Section 4 we shall deduce Theorem \ref{mink}
from the convexity and concavity properties of $\Psi_{p,q}$, again
using the same strategy that was employed in \cite{CL} in the case
$0<p \leq 1,\ q=1$. Finally, in Section 5, we discuss relations among
the various convexity and concavity results discussed here.  Roughly,
it turns out that all of the convexity statements are equivalent to one
another: There is a chain of deduction starting from any one of them,
and leading to any other.  The same is true for concavity. 
Section 5  concludes with some brief historical comments.

\medskip
{\it Acknowledgements} We thank Frank Hansen and Mary Beth~Ruskai for
helpful comments on a draft of this paper. We also thank Quanhua Xu for
a careful reading of the draft and for pointing out a number of misprints.  

\section{Convexity and concavity properties of
$\Upsilon_{p,q}$}\label{genep}

\medskip

As indicated at the end of the introduction, the main
methodological novelty of the paper lies in this section.  The proof of
convexity of $\Upsilon_{p,q}$ divides into two cases which are $1\leq
q \leq p \leq 2$ and $1\leq p \leq 2 $ with $ q> p$.

The next lemma treats the latter case, $q >p$, which is the easiest.
\begin {lm} \label{q>p} For  $1\leq p \leq 2 $ and $ q> p$, \,
$\Upsilon_{p,q} $ is convex on $\hnp$.  \end{lm}

\noindent {\bf Proof:} Since $r:= q/p \geq 1$ and since $B^*A^pB \geq
0$, we can write 
\begin{equation} \label{convex} 
\| B^*A^pB \|_r =
\sup_{\scriptsize{\begin{array}{c} \|Y\|_{r' }\leq 1, \\ Y\geq 0 \\
\end{array}} } \tr (B^*A^pB Y)\ 
\end{equation} 
where $1/r + 1/r' =1$.  Since $A^p$ is well known to be operator
convex in $A$ for $1\leq p\leq 2$, so is $B^*A^pB$.  Since the right
side of (\ref{convex}) is the supremum of a family of convex functions
(note that $Y\geq 0$ is needed here) we conclude that $\| B^*A^pB \|_r
$ is convex.  (This quantity is the $\frac{1}{r}$th power of
$\Upsilon_{p,q}(A)$, which is homogeneous of degree $r\geq 1$; therefore,
by the Remark \ref{r1} after Theorem \ref{main}, $\Upsilon_{p,q}(A) $ is
convex.)  \lanbox 

The case $q<p$ cannot be treated by such elementary means, and the
next three lemmas present the tools that we shall use. The first is a
variational formula for $p$th roots.

For all $r>1$, and all $c,x>0$, the arithmetic--geometric mean
inequality says that $$\frac{1}{r} c^r + \frac{r-1}{r}x^r
\ge cx^{r-1}\ ,$$ and so \begin{equation}\label{argeop} c =
\frac{1}{r}\inf\left\{\frac{c^r}{x^{r-1}} +(r-1)x\ :\ x>0\ \right\}\ .
\end{equation}

One may then easily verify the corresponding formula
for $0<r<1$: \begin{equation}\label{argeopcav} c =
\frac{1}{r}\sup\left\{\frac{c^r}{x^{r-1}} +(r-1)x\ :\ x>0\ \right\}\ .
\end{equation}

We shall develop this into a variational formula for
$\Upsilon_{p,q}$.  For what follows, it is first useful to
note that since $B^*A^pB$ and $A^{p/2}BB^*A^{p/2}$ have the same
spectrum, 
\begin{equation}\label{altform} 
\Upsilon_{p,q}(A) = {\rm
Tr}\left[(A^{p/2}BB^*A^{p/2})^{q/p}\right]\ .  
\end{equation}

\medskip

\begin{lm}\label{lem1} For any positive $n\times n$
matrix $A$, and with $r=p/q >1$, we have, for any $p$
\begin{equation}\label{varformp} 
\Upsilon_{p,q}(A) =
\frac{1}{r}\inf\left\{{\rm Tr}\left[A^{p/2}B\frac{1}{X^{r-1}}B^*
A^{p/2} +(r-1)X\right]\ :\ X>0\ \right\} 
\end{equation} 
where the infimum is taken over all positive $n\times n$ matrices $X$.
Likewise, for $r=p/q <1$, 
\begin{equation}\label{varformp2}
  \Upsilon_{p,q}(A) = \frac{1}{r}\sup\left\{{\rm Tr}
    \left[A^{p/2}B\frac{1}{X^{r-1}}B^* A^{p/2} +(r-1)X\right]\ :\ X>0\
  \right\} 
\end{equation} 
where the supremum is taken over all
positive $n\times n$ matrices $X$.  \end{lm} 
\medskip
\noindent{\bf Proof:} Let $C =B^*A^{p/2}$.  By continuity we may
assume that $C^*C$ is strictly positive. Then, for $ r>1$, there is a
minimizing $X$.  Let $Y= X^{1-r}$ and note that minimizing
(\ref{varformp}) with respect to $X$ is the same as minimizing
$\left\{{\rm Tr}\left(CC^* Y +(r-1) Y^{-1/(r-1)}\right)\right\}$ with
respect to $Y$.  Since the minimizing $Y$ is strictly positive, we may
replace the minimizing $Y$ by $Y+ tD$, with $D$ self adjoint, and set
the derivative with respect to $t$ equal to $0$ at $t=0$. This leads
to $\tr D[CC^* - Y^{-r/(r-1)} ] =0$. Therefore $Y^{-r/(r-1)} = CC^*$
and we are done.  The variational formula (\ref{varformp2}) is proved
in the same manner. \lanbox

\medskip

While these variational formulas may not appear to be very promising
for our purposes since, in general, the infimum of a family of convex
functions is not convex, they will  be useful on account of the following
lemma.  \begin{lm}\label{rock} If $f(x,y)$ is jointly convex, then
 $g(x)$ defined by $$g(x) = \inf_y f(x,y)$$ is convex. The analogous
 statement with convex replaced by concave and infimum replaced by
 supremum is also true.  \end{lm}

 \medskip This lemma is well known in the theory of convex functions (see
 Theorem 1 in \cite{R}), but we recall the simple proof for completeness:
 \medskip

 \noindent {\bf Proof:} For any $x_0$ and $x_1$, and any $0<\lambda<1$,
 pick $\varepsilon>0$, and choose $y_0$ and $y_1$ so that
$$f(x_0,y_0) < g(x_0)+\varepsilon \qquad{\rm and}\qquad  f(x_1,y_1) <
g(x_1)+\varepsilon\ .$$ Then \begin{eqnarray} \qquad\qquad \qquad g(
(1-\lambda)x_0+\lambda x_1) &\le& f((1-\lambda)x_0+\lambda x_1,
(1-\lambda)y_0+\lambda y_1)\nonumber\\ &\le& (1-\lambda)f(x_0,y_0)+\lambda
f(x_1,y_1)   \nonumber \\ &\le&  (1-\lambda)g(x_0) + \lambda g(x_1) +
2\varepsilon\ . \qquad\qquad\qquad \qquad\qquad\qquad \lanbox \nonumber
\end{eqnarray}

 A special case of this lemma concerns the convex function
on $\R^n \times \R^n$ given by $f(x,y)= h(x-y) +g(y)$ in which $h $ and
$g$ are convex on $\R^n$. Then $\inf_y f(x,y)$ is called the {\it infimal
convolution} of $h$ and $g$.  Its convexity accounts for the following
physical fact: When one combines two physical systems and allows heat
to flow between them so that the total energy $x $ is conserved, the
energy distributes itself so that the total entropy is maximized. Thus,
the total entropy, given by $$ S(x) = \sup_y \{ S_1(x-y) +S_2(y) \}, $$
is again a concave function of the total energy $x$, as any legitimate
entropy must be.

On account of Lemmas \ref{lem1} and \ref{rock} we can prove the
stated convexity and concavity properties of $\Upsilon_{p,q}$ if
we can prove that \begin{equation}\label{joint} (A,X) \mapsto {\rm
Tr}\left(A^{p/2}B^*\frac{1}{X^{r-1}}BA^{p/2}\right) \end{equation}
jointly convex for $1\le r \le p \le 2$, and jointly concave for $0 <
p \le r \le 1$.  In fact, for $p=r =2$  it is known \cite{LR2} that
$$(C,X) \mapsto   C^*X^{-1}C$$ is even operator convex, which gives us
more than we need in this case.

To handle the general case, we need the following joint convexity result
of Ando \cite{A}, and a concavity result in \cite{L}:

\medskip \begin{lm}\label{al} The map $$(A,B) \mapsto A^p\otimes B^{1-r}$$
on $\hnp \times \hnp$, is jointly convex for all $1 \leq r \le p \le 2$
(Ando's convexity theorem) and is jointly concave for all $0 \leq p \leq r \leq 1$
(Lieb's concavity theorem).  \end{lm}

\remark \label{requiv} 
It is well known that the convexity/concavity stated in the lemma is equivalent to
the convexity/concavity of $(A,B) \mapsto \tr A^pK^* B^{1-r}K$ for
{\it every} matrix $K \in \mn$. The difference between the two formulations is 
merely notational.
A matrix $K_{i,j}$ can equally well be regarded as a vector $K^{\rm vec}$ in
$\C^n \otimes \C^n$. Then $\tr A^pK^* B^{1-r}K =
\langle K^{\rm vec} ,\, A^p\otimes B^{1-r} \, K^{\rm vec} \rangle $.

\medskip Proofs of both results can be found in \cite{A}; see Corollary
6.2 for the concavity, and Corollary 6.3 for the convexity.  Next, we
show that Lemma \ref{al} provides the desired convexity and concavity
properties of the map defined in (\ref{joint}).

\medskip \begin{lm}\label{jcc} The map  on $\hnp \times \hnp$ defined in
(\ref{joint}) 
 is jointly convex  for all $1 \leq r \le p \le 2$ and is jointly concave
 for all $0 < p < r < 1$.
\end{lm}

\medskip \noindent{\bf Proof:} We first rewrite the right hand
side of (\ref{joint}) in a more convenient form: Define $$Z =
\left[\begin{array}{ccc} A & 0 \\ 0 & X \\ \end{array}\right]
\qquad{\rm and}\qquad K =   \left[\begin{array}{ccc} 0 & 0 \\ B & 0 \\
\end{array}\right]\ .$$ Then $$K^*Z^{1-p}K = \left[\begin{array}{ccc}
B^*X^{1-p}B & 0 \\ 0 & 0 \\ \end{array}\right]\ ,$$ and so, using
the cyclicity of the trace, $${\rm Tr}(Z^pK^*Z^{1-r}K) =   {\rm
Tr}\left(A^{p/2}B^* \frac{1}{X^{r-1}}BA^{p/2}\right)\ .$$ Note that
convexity/concavity of the left hand side in $Z$ is the same as
convexity/concavity of the right hand side in $(A,X)$.

The relation between the left hand side  and $Z^p\otimes Z^{1-r}$
is explicated in Remark \ref{requiv}. 
\lanbox

\medskip

We are now ready to prove the main result of this section:

\medskip \begin{lm}\label{upprop} For $1 \le q \le p \le 2$,
$\Upsilon_{p,q}$  is convex on $\hnp$,
 while for $0 \le p \le q \le 1$,  $\Upsilon_{p,q}$ is concave on $\hnp$.
 \end{lm}

 \noindent{\bf Proof:} By Lemma \ref{jcc}, the mapping in (\ref{joint})
 is jointly convex for $1 \le r \le p \le 2$. Then taking $r= p/q$, we
 have from Lemma \ref{lem1} and Lemma \ref{rock} $$\Upsilon_{p,q}(A) =
 \inf_X f(A,X)$$ where $f(A,X)$ is jointly convex in $A$ and $X$. The
 convexity of $\Upsilon_{p,q}$ now follows by Lemma \ref{rock}.

 The concavity statement is proved in the same way. \lanbox

 \medskip

\section{Convexity and concavity properties of $\Phi_{p,q}$ and
$\Psi_{p,q}$}\label{others}

\medskip \begin{lm}\label{phiprop}
  For $1 \le p \le 2$  and all $q\geq 1$, \ $\Phi_{p,q}$ is jointly
  convex on $(\hnp)^m$, while for $0 \le p \le q \le 1$, $\Phi_{p,q}$
  is jointly concave on $(\hnp)^m$.
 \end{lm}

 \medskip
\noindent{\bf Proof:} First, consider the case $m=2$ and define 
$${\cal A}
=  \left[\begin{array}{ccc} A_1&0\\
 0&A_2\\ \end{array}\right] \qquad{\rm and}\qquad
\sigma =  \left[\begin{array}{ccc} 0&I\\ I&0\\ \end{array}\right] \ \ .
$$
Then 
$${\cal A}^p + \sigma {\cal A}^p\sigma =
 \left[\begin{array}{ccc} A_1^p+A_2^p&0\\ 0&A_1^p+A_2^p\\
 \end{array}\right]
$$
But 
$${\cal A}^p + \sigma {\cal A}^p\sigma = 2\left({I+\sigma\over
2}\right){\cal A}^p\left({I+\sigma\over 2}\right) + 2\left({I-\sigma\over
2}\right){\cal A}^p\left({I-\sigma\over 2}\right)
$$ 
Now define 
$$\Pi_\pm
= {I\pm \sigma\over 2}
$$ 
and observe that these are complementary
orthogonal projections.  Therefore, the $q/p$th power of the sum
on the right is simply the sum of the $q/p$th powers, and hence
\begin{equation}\label{split} 2\tr\left[(A_1^p + A_2^p)^{q/p}\right]
= 2^{q/p}\tr \left[\left(\Pi_+ {\cal A}^p\Pi_+\right)^{q/p}\right]
+ 2^{q/p}\tr \left[\left(\Pi_- {\cal A}^p\Pi_-\right)^{q/p}\right]
\end{equation}

By Lemmas \ref{q>p} and \ref{upprop}, applied with $A = {\cal A}$
and $B = \Pi_\pm$, we see that each of the two terms on the right hand
side of (\ref{split}) is a convex function of ${\cal A}$, for $1 \le
p \le 2\ , \ q\geq 1$, and concave for $0< p \le q \le 1$.  Of course,
convexity/concavity in ${\cal A}$ amounts to joint convexity/concavity
in $A_1$ and $A_2$.

Apart from a factor of $2$, the left hand side of (\ref{split}) is
the $q$th power of $\Phi_{p,q}(A_1,A_2)$.  Again, using Remark
\ref{r1}, and taking the $q$th root of the right hand side, we preserve
the convexity and concavity properties since the right hand side is
homogeneous of degree $q$.  Hence $\Phi_{p,q}(A_1,A_2)$ has the stated
convexity and concavity properties.  This concludes the proof for $m=2$.

One can easily iterate this procedure to obtain the result for all dyadic
powers $m = 2^k$, and hence for all $m$.  \lanbox

\medskip

We now turn to $\Psi_{p,q}$:

\medskip

\begin{lm}\label{psiprop} For $1 \le p \le 2$ and $q\geq 1$,  $\Psi_{p,q}$
is  convex on $\hmnp$,
 while for $0 \le p \le q \le 1$,  $\Psi_{p,q}$ is  concave on $\hmnp$.
\end{lm}
\medskip

The proof turns on Lemma \ref{phiprop} and a method for writing
partial traces as averages, which originates with Uhlman \cite{U}. To
recall this, fix some orthonormal basis $\{e_1,e_2,\dots,e_n\}$ of
$\C^n$. Let ${\cal G}$ be the group of unitary operators $W$ such
that for some permutation $\pi$ of $\{1,\dots,n\}$, and some function
$s:{1,2,\dots,N}\mapsto \{0,1\}$, \begin{equation}\label{wdef} We_j =
(-1)^{s(j)}e_{\pi(j)}\qquad{\rm for\ all}\ j\ .  \end{equation}

It is clear that these form a group, and that the cardinality of the group
is $2^n n!$.  The point of the definition is that any operator on $\C^n$
that commutes with every element of ${\cal G}$ is necessarily a multiple
of the identity, and so,  for any operator $A$ on $\C^n\otimes \C^n$,
\begin{equation}\label{uhl}
 \frac{1}{n}\left[{\rm Tr}_2(A)\right]\otimes I_{n\times n} =  \frac{1}{2^n n!}
\sum_{W\in{\cal G}} [I\otimes W^*]A \, [I\otimes W]\end{equation} \bigskip

\noindent{\bf Proof of Lemma \ref{psiprop}:}\quad
Applying \ref{uhl}, we have 
\begin{eqnarray} 
&&{\rm
Tr}_1\left[\left({\rm Tr}_2A^p\right)^{q/p}\right] = n^{q/p-1}
{\rm Tr}_{1}\tr_{2}\left[\left(\frac{1}{n}[{\rm Tr}_2A^p]\otimes I_{n\times
n}\right)^{q/p} \right] \nonumber\\ 
&&=n^{q/p-1}{\rm Tr}_{1}\tr_{2}\left[\left(
\frac{1}{2^n n!}\sum_{W\in{\cal G}} [I\otimes W^*]A^p \, [I\otimes W]
\right)^{q/p}\right] \nonumber\\ 
&&=n^{q/p-1}\left(\frac{1} {2^n
n!}\right)^{q/p}{\rm Tr}_{1}\tr_{2}\left[\left( \sum_{W\in{\cal G}}\left(
[I\otimes W^*]A\, [I\otimes W] \right)^p \right)^{q/p}\right] \ .\nonumber
\end{eqnarray}
Upon taking the $q$th root, the  result  follows directly from Lemma
\ref{phiprop}.\quad \lanbox \bigskip

\remark While the proof shows that the conclusion of Lemma
\ref{psiprop} follows from Lemma \ref{phiprop}, the reverse implication is
even more elementary: To see this, suppose that the matrix $A$ in Lemma
\ref{psiprop} is the block diagonal matrix whose $j$th diagonal block is
$A_j$.  Then, clearly, 
$$\Psi_{p,q}(A) = \Phi_{p.q}(A_1,A_2,\dots,A_m)\ .
$$

We next show that for $p>2$, $\Phi_{p,q}$ is neither concave nor convex,
not even separately in each $A_j$. The argument is a simple adaptation
of the one we gave in \cite{CL} for $q=1$. However, the argument has
other uses, and so it will be useful to give it here explicitly.

\medskip

\begin{lm}\label{noncon} For $p>2$, $\Phi_{p,q}$ is neither convex nor
  concave, even separately,
 for any values of $q\ne p$.
\end{lm}

\medskip \noindent{\bf Proof:}  The basis of the proof is a simple
Taylor expansion. By simple differentiation one finds that for any $A,
B\in \hnp$, \begin{equation}\label{taylor} \Phi_{p,q}(tA,B) = \|B\|_q
+\frac{t^p}{p}   \|B\|_q^{1-q} \, \tr A^pB^{q-p}  + O(t^{2p})\ .
\end{equation}

Keeping $B$ fixed, but replacing $A$ by $A_1$, $A_2$ and $(A_1+A_2)/2$,
we find \begin{eqnarray}\label{expand} &&\frac{1}{2} \Phi_{p,1}(tA_1,B) +
\frac{1}{2} \Phi_{p,1}(tA_2,B) -\Phi_{p,1}\left(t\frac{A_1+A_2}{2},\,
B \right)  = \nonumber\\ &&\frac{t^p}{p}\, \|B\|_q^{1-q}
\left[\frac{1}{2} \tr \left(A_1^p  B^{q-p} \right)+ \frac{1}{2} \tr
\left(A_2^p B^{q-p}\right)  - \tr \left(\left(\frac{A_1+A_2}{2}\right)^p
\,    B^{q-p} \right) \right] + O(t^{2p})\ .\nonumber 
\end{eqnarray}

Now if $p>2$, $A\mapsto A^p$ is not operator convex, and so we can find
$A_1$ and $A_2$ in $\hnp$ and a unit vector $v$ in $\C^n$ such that
\begin{equation}\label{rone} 
\frac{1}{2}\langle v, A_1^p v \rangle  +
\frac{1}{2}\langle v, A_2^p v \rangle - \left\langle v, \left(\frac{A_1 +
A_2}{2}\right)^p v \right\rangle \end{equation} is strictly negative.
However, since $A\mapsto \tr A^p$ is convex, there is necessarily
another unit vector
$w$, so that if we replace $v$ by $w$, this quantity becomes positive.

Hence  for $q>p$, take $B$ to be the rank one projection onto $v$. Then the
bracketed quantity on the right hand side of (\ref{expand}) is nothing
other than (\ref{rone}), and so the left had side of (\ref{expand})
is strictly negative. This shows that $\Phi_{p,q}$ cannot be convex for
such $p$ and $q$.  Replacing $v$ by $w$, the sign changes, and we see
that it cannot be concave either.

For $q<p$, take $B = (\varepsilon I + vv^*)^{-1}$ for some  small value of
$\varepsilon$, so that $B^{q-p}$ is essentially the orthogonal projection
onto $v$, and then argue as before. \lanbox

\medskip We have now completed all the steps needed for the proof of
the main theorem:

\medskip

\noindent{\bf Proof of Theorem \ref{main}:} Lemmas \ref{q>p},
\ref{upprop}, \ref{phiprop} and \ref{psiprop} respectively establish
the convexity of $\Upsilon_{p,q}$, $\Phi_{p,q}$ and $\Psi_{p,q}$ for
$1\le p\le 2$,  $q\geq 1$,  and the concavity of these functionals
for $0 < p \le q \le 1$.  Since we deduced the convexity/concavity
properties of $\Phi_{p,q}$ from those of $\Upsilon_{p,q}$ in the
proof of Lemma \ref{phiprop}, it follows that if $\Phi_{p,q}$ is not
convex or concave for some values of $p$ and $q$, then neither is
$\Upsilon_{p,q}$. Likewise, in the remark following the proof of Lemma
\ref{psiprop}, we have observed that $\Phi_{p,q}$ and $\Psi_{p,q}$
have the same convexity/concavity properties. Therefore, since Lemma
\ref{noncon} says that $\Phi_{p,q}$ is neither convex nor concave for
$p>2$,  $q\ne p$, the same is true of $\Upsilon_{p,q}$ and $\Psi_{p,q}$ \lanbox

\section{Minkowski's inequality for trace norms on a product of three
  Hilbert spaces}\label{others}

\medskip

\noindent{\bf Proof of Theorem \ref{mink}:}\quad First we reduce to
the case $q=1$: Suppose we have shown that 
\begin{equation}\label{q=1}
{\rm Tr}_3\left[{\rm  Tr}_2({\rm Tr}_1 A)^{p}\right]^{1/p} \le {\rm
Tr}_3\tr_1 \left[({\rm Tr}_2 A^{p})^{1/p} \right]  \ 
\end{equation}
for all positive $A$, and all $1 \le p \le 2$.  Then replacing $A$ by
$A^q$, and $p$ by $r = p/q$, we obtain (\ref{minkpq}). The same works
for the reverse inequality for $0 \le p \le 1$ except now we require
$q\ge p$  to have $0 \le r \le 1$.  Therefore, it remains to prove
(\ref{q=1}). Consider first the case $1 \le p \le 2$.

Suppose that the dimension of $\h_1$ is $n$. The left hand side of
(\ref{q=1}) can be written in terms of $\Psi_{p,1}$: 
\begin{eqnarray}
{\rm Tr}_3\left[{\rm  Tr}_2({\rm Tr}_1A)^{p}\right]^{1/p} &=& {\rm
Tr}_{3}\tr_{1}\left[{\rm Tr}_2\left( {1\over n} I_{\h_1}\otimes {\rm Tr}_1
A\right)^{p}\right]^{1/p} \nonumber\\ 
&=& \Psi_{p,1}\left({1\over
n} I_{\h_1}\otimes {\rm Tr}_1 A\right)\nonumber
\end{eqnarray}
where the pair of spaces in the definition of $\Psi_p$ is taken to be
$\h_1\otimes\h_3$ and $\h_2$.  In the following, it will be useful to
let $I_2$ denote $I_{\h_2}$, to let $I_{1,3}$ denote $I_{\h_1\otimes \h_3}$,
and so forth.

Then by (\ref{uhl}) and the convexity of $\Psi_{p,q}$ for $1\le q\le
p \le 2$, 
\begin{eqnarray}\label{usepsi} \Psi_{p,1}\left({1\over
n} I_{\h_1}\otimes {\rm Tr}_1 A\right) 
&=& \Psi_{p,1}\left(
\frac{1}{2^n n!}\sum_{W\in{\cal G}} [W^*\otimes I_{2,3}]A\,
[W\otimes I_{2,3})  \right] \nonumber\\ 
&\le&\frac{1}{2^n
n!} \sum_{W\in{\cal G}}\Psi_{p,1}\left([W^*\otimes I_{2,3}]A \,
[W\otimes I_{2,3}] \right)\nonumber
\end{eqnarray} 
Now, by
the definition of $\Psi_{p,1}$, and the fact that $W$ is unitary,
$$
\Psi_{p,1}\left([W^*\otimes I_{2,3}]A\, [W\otimes I_{2,3}] \right) =
{\rm Tr}_{1,3}\left[{\rm Tr}_2\left([W^*\otimes I_{2,3}]A^p \, [W\otimes
I_{2,3}]\right) \right]^{1/p} \ .
$$ 
However, 
$$
{\rm Tr}_2\left([W^*\otimes
I_{2,3}]A^p\, [W\otimes I_{2,3}]\right)  = [W^*\otimes I_3] {\rm Tr}_2(A^p)\,
[W\otimes I_3]\ ,
$$ 
and hence, again by the unitarity of $W$, 
$$
\left[{\rm
Tr}_2\left([W^*\otimes I_{2,3}]A^p \, [W\otimes I_{2,3}]\right)\right]^{1/p}
= [W^*\otimes I_3]\left[{\rm Tr}_2(A^p)\, \right]^{1/p} [W\otimes
I_3]\ .
$$ 
After taking the trace over $\h_1 \otimes \h_3$, we obtain 
$$
\Psi_{p,1}\left([W^*\otimes I_{2,3}]A \, [W\otimes I_{2,3}] \right) =
{\rm Tr}_3{\rm Tr}_1\left({\rm Tr}_2
    A^p\right)^{1/p}\ .
$$ 
Thus  each term in the average in (\ref{usepsi}) is the same, and
we obtain (\ref{minkpq}). Note that since $\Psi_{p,1}$ is concave in
the case $0< p \le q \le 1$, the inequality in (\ref{usepsi}) reverses,
and with it the inequality in (\ref{minkpq}). \lanbox


\section{What follows from what}

In the chain of deduction leading up to the proof of Theorem \ref{main},
our starting point was Ando's convexity theorem   and  Lieb's concavity
theorem quoted in Lemma \ref{al} We then deduced convexity/concavity
properties of $\Upsilon_{p,q}$ from this, and then convexity/concavity
properties of  $\Phi_{p,q}$ from those of $\Upsilon_{p,q}$. We have
already noted that the convexity/concavity properties of $\Phi_{p,q}$
and $\Psi_{p,q}$ are identical.

It is interesting to observe that this chain of deduction can be brought
full circle: Given the convexity of $\Phi_{p,1}$ for $1\le p\le 2$, and
the concavity for $0< p \le 1$, one can prove Ando's convexity theorem
and Lieb's concavity theorem.  This can be done using some ideas of
Bekjan.  In \cite{B}, Bekjan proved the joint concavity of $\Phi_{p,1}$
for  $-1 < p < 0$. He did this by first adapting Epstein's proof of his
inequality to negative values of $p$, and then applying the method used in
\cite{CL} to prove the joint concavity of $\Phi_{p,1}$ for $0< p \le 1$.

His main interest in this result seems to have been the consequence
that he derived from it by using the Taylor expansion (\ref{taylor})
at $q=1$, where  it simplifies to \begin{equation} \Phi_{p,1}(tA,B)
= \tr B +\frac{t^p}{p} \tr B^{1-p}A^p  + O(t^{2p})\ , \end{equation}
for $A,B \in \hnp$.
As Bekjan showed, since the constant term is affine, one can deduce from this formula and the 
joint concavity  of $(A,B) \mapsto   \Phi_{p,1}(A,\,B) $ for $-1 \le p \le 0$  
the joint convexity on $\hnp\times \hnp$ of $(A,B) \mapsto   \tr A^p\, B^{1-1/p}$ for $1\le p \le 2$. 
Now that we know that the joint convexity of   $(A,B) \mapsto   \Phi_{p,1}(A,\,B) $ for $1 \le p \le 2$,
the argument can be made even more direct:

\begin{lm} \label{bekjan}  Consider the maps $(A,B) \mapsto   \Phi_{p,1}(A,\,B) $  and $(A,B) \mapsto   \tr A^p\, B^{1-1/p}$
on $\hnp \times \hnp$.   
For $ 0\leq p \leq q \leq 1$ the joint
  concavity of $(A,B) \mapsto \Phi_{p,1}(A,\,B) $ implies the joint concavity of $(A,B) \mapsto \tr
  A^p\, B^{1-1/p}$.  For $1\leq q\leq p \leq 2$ the joint convexity
  of $(A,B) \mapsto \Phi_{p,1}(A,\,B)$ implies the joint convexity of $(A,B) \mapsto \tr A^p\, B^{1-1/p}$.
\end{lm}

\noindent {\bf Proof:}  It suffices to prove midpoint
convexity/concavity. Fix $A_1,A_2$ and $B_1,B_2$ in $\hnp$,
and for $t>0$ define \begin{equation} \label{fdef} f(t) =
\frac{1}{2} \Phi_{p,1}(tA_1,B_1) + \frac{1}{2} \Phi_{p,1}(tA_2,B_2)
-\Phi_{p,1}\left(t\frac{A_1+A_2}{2},\,   \frac{B_1+B_2}{2} \right)
\end{equation} Note that $f(0) =0 $ and that by Theorem \ref{main}, $f(t)
\geq 0$ (resp. $\leq 0$) for $1\leq p\leq 2$ (resp. $0\leq p \leq 1$).

We now apply the Taylor expansion (\ref{taylor}) to approximate each of
the three terms in the definition of $f$. Since the constant term in
the Taylor expansion is affine, it drops out (this is the reason for the
restriction to $q=1$ here), and consequently, to leading order in $t$,
\begin{equation}
  f(t) \sim    \frac{t^p}{p}   \left(\frac{1}{2} \tr \left[A_1^p B_1^{1-p}
  \right]+
\frac{1}{2} \tr \left[A_2^p B_2^{1-p}\right]  - \tr
\left[\left(\frac{A_1+A_2}{2}\right)^p \,    \left(\frac{B_1+B_2}{2}
\right)^{1-p} \right]\right)\end{equation} The   quantity on
the right that multiplies $t^p/p$ 
must have the same sign as $f(t) $ for small $t$, which proves
the convexity/concavity.  \lanbox \medskip

This leads to the conclusion that $(A,B) \mapsto \tr A^p B^{1-p}$ is
jointly convex/concave on $\hnp\times \hnp$ depending on the value of $p$.  Bekjan further
observed that one can ``put the $q$ back'' into this inequality,
thereby obtaining (with $r=p/q$) the more general result that $(A,B)
\mapsto \tr A^p B^{1-r}$ is convex for $1 \le r \le p \le 2$.  Here is
Bekjan's argument: Observe that $0 \le (1-r)/(1-p) \le 1$, so that
$A \mapsto A^{(1-r)/(1-p)}$ is concave. Also, $-1\le 1-p \le 0$,   
so $A \mapsto A^{1-p}$ is monotone non-increasing on $\hnp$. Thus, for any
$B_1,B_2\in \hnp$, and any $\lambda\in (0,1)$, 
\begin{eqnarray}
(\lambda B_1 + (1-\lambda)B_2)^{1-r} 
&=& \left( [\lambda B_1 +
(1-\lambda)B_2]^{(1-r)/(1-p)}\right) ^{1-p}\nonumber
\\ &\le& \left(
\lambda B_1 ^{(1-r)/(1-p)} + (1-\lambda)B_2^{(1-r)/(1-p)}\right)^{1-p}\
.\nonumber
\end{eqnarray} 
Therefore, 
\begin{multline} 
\tr\left[\left(\lambda A_1 + (1-\lambda)A_2\right)^p \left(\lambda B_1 +
(1-\lambda)B_2\right)^{1-r} \right] \\
\le \tr \left[
\left(\lambda A_1 + (1-\lambda)A_2\right)^p      \left(\lambda B_1
^{(1-r)/(1-p)} + (1-\lambda)B_2^{(1-r)/(1-p)}\right)^{1-p}\right]\ .
\end{multline} 
Now applying the convexity of
$(A,B)\mapsto \tr A^pB^{1-p}$, Bekjan obtains 
\begin{multline} \qquad\qquad
\tr \left[\left(\lambda A_1 + (1-\lambda)A_2\right)^p \left(\lambda B_1 +
(1-\lambda)B_2\right)^{1-r} \right]  \\
\leq \lambda
\tr \left[A_1^p B_1^{1-r}\right]+ (1-\lambda) \tr \left[A_2^p
B_2^{1-r}\right]\ , \qquad\qquad
\end{multline} 
which is the desired
convexity conclusion.  In the same way, one extends the concavity result
to $0 \le p \le r \le 1$.  Bekjan was aware that the concavity result
one obtains in this way is equivalent to Lieb's concavity theorem,
but he was either unaware of Ando's convexity theorem, or unaware that
his trace convexity theorem is actually equivalent to Ando's result.
The equivalence is given by the following lemma:

\begin{lm}  \label{equiv} For any real numbers $s,t$ these  two statements
are equivalent: \\ 
  (1.)  \ $(A, \,  B) \mapsto  \tr \left( A^s\, B^t\right)$  is
  concave/convex on $\hnp\times \hnp$. \\ (2.) \ $(A,B)\mapsto A^s\otimes B^t$ is operator
  concave/convex $\hnp\times \hnp$.
\end{lm}

\noindent {\bf Proof:}  {\em (2.)} easily implies {\em (1.)}  because, as stated in 
Remark \ref{requiv},
{\em (2.)} is equivalent to the 
concavity/convexity of
$(A, \,  B) \mapsto  \tr \left( A^s K^*\, B^t K \right)$ for every matrix $K$. To obtain
{\em (1.)} simply take $K= I$.     

To derive {\em (2.)} from {\em (1.)},  it suffices to show that {\em (1.)} implies that
for each bounded operator $K$ on $\h$, $(A, \, B) \mapsto \tr \left(
  A^s K B^t K^* \right)$ is concave/convex, as explained in section
\ref{genep}. If $K$ is a unitary, $U$, this follows immediately from
{\em (1.)} since $(UBU^*)^t = UB^tU^*$ and concavity/convexity in $(A,B)$ is
equivalent to concavity/convexity in $(A,UBU^*)$. To reduce the general
case to the unitary case we may assume, without loss, that $K$ is a
contraction.  Then $K$ can be dilated to a unitary on the larger space
$\h\oplus \h$. That is, there is an operator $L$ such that ${\cal W} =
\left[\begin{array}{ccc}\phantom{-}K & L \\ -L & K\\ \end{array} \right]
$ is unitary, namely, $L=U(I-|K|^2)^{1/2}$, where $K=U|K|$ is the polar
decomposition of $K$. Then, with ${\cal A} = \left[\begin{array}{ccc}A
    & 0\\ 0 & 0 \\ \end{array} \right] $ and ${\cal B} =
\left[\begin{array}{ccc}B& 0\\ 0 & 0 \\ \end{array} \right] $
\begin{equation}
 \tr_\h \left( A^s\, K \, B^t\,  K^*\right) =
\tr_{\h \oplus \h} \left( {\cal A}^s\, {\cal W}\, {\cal B}^t \, {\cal
W}^*\right) , \end{equation} and the general case follows from the
unitary case. \lanbox \medskip

Thus, the joint convexity of $(A,B) \mapsto \tr (A^p\, B^{1-r})$ for $1
\le r \le p \le 2$, proved in \cite{B} is actually equivalent to Ando's
convexity theorem, quoted in Lemma \ref{al}.
 This
equivalence is significant for us because it brings our chain of
implications full circle, back to the convexity and concavity theorems for
tensor products. In Lemma \ref{upprop} we deduced the convexity/concavity
properties of $\Upsilon_{p,q}$ from those of the tensor products. In
Lemma \ref{phiprop} we deduced the convexity/concavity properties of
$\Phi_{p,q}$ from $\Upsilon_{p,q}$.  In Lemma \ref{psiprop} and the
remark following it we proved that $\Phi_{p,q}$ and $\Psi_{p,q}$ have
the same convexity/concavity properties.

The chain of implications   is summarized in the following diagram,
in which the labels on the arrows give the relevant  lemmas.  \bigskip
\[ \begin{CD} 
A^p\otimes B^{1-p/q} @>{\rm Lemma}\ \ref{upprop} >>
  \Upsilon_{p,q}\\ @A{\rm
    Lemma}\ \ref{equiv}AA             @VV{\rm Lemma}\ \ref{phiprop} V\\
  \tr\left(A^p \, B^{1-p/q} \right) @<<{\rm Lemma} \
  \ref{bekjan}<\Phi_{p,q} @=\Psi_{p,q} \ {\scriptsize{\text { (Remark
        \ after\ Lemma 3.2)} } }\\ 
\end{CD} \]


\subsection{ Historical Remarks:} We end this paper with some historical
comments on the large literature on related matrix convexity and
concavity theorems. This is by no means meant to be a complete survey.

Many lines of investigation in this field can be traced back to
the Wigner--Yanase papers \cite{WY1,WY2} that drew attention to the
concavity of $A\mapsto \tr (A^p K A^{1-p} K^*)$ and proved it for $p=1/2$.
Their motivation was that in a quantum system, some states are easier to
measure than others; if a density matrix $\rho$ commutes with a conserved
quantity (say the energy) then it is easy to measure, and otherwise
not. Thus, while the von Neuman entropy of any pure state $\rho$ is
zero, some pure states have a higher information content than others
-- namely those that are not functions of the conserved quantities. In
the presence of one conserved quantity, represented by the self-adjoint
operator $K$, the {\it  Wigner--Yanase skew information} of $\rho$ is
$$-\tr \left([\sqrt{\rho} ,K]^2\right)\ .$$ Note that this is $$I(\rho)
= \tr K^2\rho - \tr \sqrt{\rho}K \sqrt{\rho}K\ .$$ In \cite{WY1},
Wigner and Yanase listed a number of requirements that they considered
to be necessary for a measure of information.  They conjectured these
to hold for $I(\rho)$.  One of  them was the convexity of $I(\rho)$,
and this was proved in \cite{WY2}.  Another was a certain subadditivity
property, and this has only very recently been shown to be {\it false}
by Hansen \cite{Han}. The hope that subadditivity would hold under
additional assumptions was shown to be false by Seiringer \cite{S2}

The general case $\tr (A^p K A^{q} K^*)$ with $p,q>0, \ p+q\leq 1$ (and
$K$ not necessarily self-adjoint) was done in \cite{L} in connection with
a proof of strong subadditivity of entropy \cite{LR}.  This concavity
result was the impetus to several papers and several different proofs.
The earliest after \cite{L} was by Epstein \cite{E}
who derived the concavity, and more, by using the complex-analytic theory
of Herglotz functions. Not only did he prove the concavity in \cite{L}, he
also proved the concavity of $\Upsilon_{p,1}$ for $0\leq p \leq 1$. (Paper
\cite{E} also led to other results not directly related to our line of
inquiry here.)  Another one is Uhlmann's fundamental result \cite{U2}
on the monotonicity of the relative entropy in a von Neuman algebra,
building on Araki's extension \cite{Ar} of the notion of relative
entropy to this more general setting. Another useful extension of strong
subadditivity and related concavity inequalities is in \cite{LS}; this
was applied in \cite{S}. .

Ando's work \cite{A} was also partially motivated by \cite{L} and
gave another proof of the concavity and many new concavity/convexity
results, some of which are used here.  His was the first proof using
purely real variable methods.  
He realized that the fact that the
{\it geometric mean} of two positive matrices $A$ and $B$, namely
$A^{1/2}(A^{-1/2}BA^{-1/2})^{1/2}A^{1/2}$, is a jointly concave
function of $A$ and $B$, as is the {\it harmonic mean}, namely
$2(A^{-1}+B^{-1})^{-1}$, is intimately related to the
convexity/concavity properties of $(A,B) \mapsto A^p\otimes B^{1-r}$,
or, equivalently, of $A\mapsto \tr (A^p K A^{1-r} K^*)$.

This result in \cite{L} has been generalized in many ways, e.g., Kosaki
\cite{K} found new proofs based on interpolation theory that provided an
analog of the concavity theorem in a general von Neuman algebra setting.

Since the original proof of strong subadditivity of quantum entropy,
there have been many others (see Ruskai's review \cite{RU1}).  One proof
was given in our paper \cite{CL}, in which strong subadditivity was
derived from a Minkowski type inequality for traces, the $0\le p \le 1$,
$q=1$ version of Theorem \ref{mink}.  Another proof \cite{D} that is
relevant for the present paper came from recent results in the
theory of operator spaces; this proof also uses a Minkowski inequality,
but in Pisier's non-commutative $L^q(L^p)$ norm (\ref{pisier}) discussed
in the introduction.  One interesting remaining issue is the condition
for equality in our convexity/concavity theorems, similar to the related
issue of equality for strong subadditivity discussed in \cite{RU1,Ha}.

\end{document}